\documentclass{amsart}
\usepackage{color}
\usepackage{amsmath} 
\usepackage{amssymb}
\usepackage{mathrsfs}

\newtheorem{theorem}{Theorem}[section] 
\newtheorem{claim}[theorem]{Claim}

\theoremstyle{definition}
\newtheorem{definition}[theorem]{Definition}

\newtheorem{problem}[theorem]{Problem}
\newtheorem{observation}[theorem]{Observation} 

\newtheorem{question}[theorem]{Question}

\theoremstyle{remark}
\newtheorem{remark}[theorem]{Remark}
\newtheorem{notation}[theorem]{Notation}

\newcommand{\rest}{{\restriction}}

\newcommand{\wilog}{{\rm without loss of generality}}

\newcommand{\then}{{\underline{then}}}
\newcommand{\when}{{\underline{when}}}
\newcommand{\Then}{{\underline{Then}}}

\newcommand{\mn}{{\medskip\noindent}}
\newcommand{\sn}{{\smallskip\noindent}}

\newcommand{\cA}{{\mathcal A}}

\newcommand{\cD}{{\mathcal D}}
\newcommand{\gd}{{\mathfrak d\/}}

\newcommand{\cF}{{\mathcal F}}

\newcommand{\bbP}{{\mathbb P}}
\newcommand{\cP}{{\mathcal P}}

\newcommand{\bbQ}{{\mathbb Q}}

\newcommand{\cT}{{\mathcal T}}

\newcount\skewfactor
\def\mathunderaccent#1#2 {\let\theaccent#1\skewfactor#2
\mathpalette\putaccentunder}
\def\putaccentunder#1#2{\oalign{$#1#2$\crcr\hidewidth
\vbox to.2ex{\hbox{$#1\skew\skewfactor\theaccent{}$}\vss}\hidewidth}}
\def\name{\mathunderaccent\tilde-3 }

\newenvironment{PROOF}[2][\proofname.]
   {\begin{proof}[#1]\renewcommand{\qedsymbol}{\textsquare$_{\rm #2}$}}
   {\end{proof}}

\begin{document}

\title {Diamonds}
\author {Saharon Shelah}
\address{Einstein Institute of Mathematics\\
Edmond J. Safra Campus, Givat Ram\\
The Hebrew University of Jerusalem\\
Jerusalem, 91904, Israel\\
 and \\
 Department of Mathematics\\
 Hill Center - Busch Campus \\ 
 Rutgers, The State University of New Jersey \\
 110 Frelinghuysen Road \\
 Piscataway, NJ 08854-8019 USA}
\email{shelah@math.huji.ac.il}
\urladdr{http://shelah.logic.at}
\thanks{Research supported by the United-States-Israel Binational
Science Foundation (Grant No. 2002323), Publication No. 922.
The author thanks Alice Leonhardt for the beautiful typing.}

\date{Oct.7, 2009}

\begin{abstract}
If $\lambda = \chi^+ = 2^\chi > \aleph_1$ then diamond on
$\lambda$ holds.  Moreover, if $\lambda = \chi^+ = 2^\chi$ and $S
\subseteq \{\delta < \lambda:\text{cf}(\delta) \ne \text{ cf}(\chi)\}$
is stationary then $\diamondsuit_S$ holds.  Earlier this was known
only under additional assumptions on $\chi$ and/or $S$.
\end{abstract}

\maketitle
\numberwithin{equation}{section}
\newpage

\section {Introduction} 

We prove in this paper several results about diamonds.  Let us
recall the basic definitions and sketch the (pretty long) history of
the related questions.

The diamond principle was formulated by Jensen, who proved that it
holds in $\bold L$ for every regular uncountable cardinal $\kappa$ and
stationary $S \subseteq \kappa$.  This is a
prediction principle, which asserts the following:
\begin{definition}
\label{0z.1}  
$\diamondsuit_S$ (the set version).

Assume $\kappa = \text{ cf}(\kappa) > \aleph_0$ 
and $S \subseteq \kappa$ is stationary,
$\diamondsuit_S$ holds \when \, there is a sequence $\langle
A_\alpha:\alpha \in S\rangle$ such that $A_\alpha \subseteq \alpha$ for
every $\alpha \in S$ and the set $\{\alpha \in S:A \cap \alpha =
A_\alpha\}$ is a stationary subset of $\kappa$ for every $A \subseteq \kappa$.
\end{definition}

The diamond sequence $\langle A_\alpha:\alpha \in S\rangle$ guesses
enough (i.e., stationarily many) initial segments of every $A \subseteq
\kappa$.  Several variants of this principle were formulated, for
example:
\begin{definition}
\label{0z.4}  
$\diamondsuit^*_S$.

Assume $\kappa = \text{ cf}(\kappa) > \aleph_0$ 
and $S$ is a stationary subset of $\kappa$.  Now
$\diamondsuit^*_S$ holds \when \, there is a sequence $\langle
{\cA}_\alpha:\alpha \in S\rangle$ 
such that each ${\cA}_\alpha$ is a subfamily
of ${\cP}(\alpha),|{\cA}_\alpha| \le |\alpha|$ and for every $A
\subseteq \kappa$ there exists a club $C \subseteq \kappa$ such that $A \cap
\alpha \in {\cA}_\alpha$ for every $\alpha \in C \cap S$.
\medskip

We know that $\diamondsuit^*_S$ holds in $\bold L$ for
every regular uncountable $\kappa$ and stationary $S \subseteq \kappa$.  Kunen
proved that
$\diamondsuit^*_S \Rightarrow \diamondsuit_S$.  Moreover, if $S_1
\subseteq S_2$ are stationary subsets of $\kappa$ then
$\diamondsuit^*_{S_2} \Rightarrow \diamondsuit^*_{S_1}$ (hence
$\diamondsuit_{S_1}$).   But the assumption
$\bold V = \bold L$ is heavy.  Trying to avoid it, we can walk in
several directions.  On weaker relatives see \cite{Sh:829} and
references there.  We can also use other 
methods, aiming to prove the diamond without assuming $\bold V = \bold L$.
\end{definition}

There is another formulation of the diamond principle, phrased via
functions (instead of sets).  Since we use this version in our proof,
we introduce the following:
\begin{definition}
\label{d.21}  
$\diamondsuit_S$ (the functional version).  

Assume $\lambda = \text{ cf}(\lambda) > \aleph_0,S \subseteq
\lambda,S$ is stationary.  $\diamondsuit_S$ holds if there exists a
diamond sequence $\langle g_\delta:\delta \in S\rangle$ which means
that $g_\delta \in {}^\delta \delta$ for every $\delta \in S$, and 
for every $g \in
{}^\lambda \lambda$ the set $\{\delta \in S:g \rest \delta =
g_\delta\}$ is a stationary subset of $\lambda$.
\end{definition}

By Gregory \cite{Gre76} and Shelah \cite{Sh:108} we know that
assuming $\lambda = \chi^+ = 2^\chi$ and $\kappa = \text{
cf}(\kappa) \ne \text{ cf}(\chi),\kappa < \lambda$, and GCH holds (or
actually just $\chi^\kappa = \chi$ or $(\forall \alpha < \chi)
(|\alpha|^\kappa < \chi) \wedge \text{ cf}(\chi) < \kappa)$, then
$\diamondsuit^*_{S^\lambda_\kappa}$ holds (recall that $S^\lambda_\kappa =
\{\delta < \lambda:\text{cf}(\delta) = \kappa\}$).

We have got also results which show that the failures of the diamond
above a strong limit cardinal are limited.  For instance, if $\lambda
= \chi^+ = 2^\chi > \mu$ and $\mu > \aleph_0$ is strong limit, then (by
\cite{Sh:460}) the set $\{\kappa < \mu:
\diamondsuit^*_{S^\lambda_\kappa}$ fails$\}$ is bounded in $\mu$
(recall that $\kappa$ is regular).
Note that the result here does not completely subsume the earlier
results when $\lambda = 2^\chi = \chi^+$ as we get ``diamond on every
stationary set $S \subseteq \lambda \backslash
S^\lambda_{\text{cf}(\chi)}$", but not $\diamondsuit^*_S$; this is
inherent as noted in \ref{2b.7}.  In 
\cite{Sh:829}, a similar, stronger result is proved for
$\diamondsuit_{S^\lambda_\kappa}$: for every $\lambda = \chi^+ = 2^\chi >
\mu,\mu$ strong limit for some \underline{finite} 
${\gd} \subseteq \text{ Reg } \cap
\mu$, for every regular $\kappa < \mu$ not from ${\gd}$ we have
$\diamondsuit_{S^\lambda_\kappa}$, and even $\diamondsuit_S$
for ``most" stationary $S \subseteq S^\lambda_\kappa$.  In fact, for
the relevant good stationary sets $S \subseteq S^\lambda_\kappa$ 
we get $\diamondsuit^*_S$.  Also
 weaker related results are proved there 
for other regular $\lambda$ (mainly $\lambda =\text{ cf}(2^\chi))$.
\bigskip

\noindent
The present work does not resolve:
\begin{problem}
\label{0z.7}  
Assume $\chi$ is singular, $\lambda = \chi^+
= 2^\chi$, do we have $\diamondsuit_{S^\lambda_{\text{cf}(\chi)}}$?
(You may even assume G.C.H.).

However, the full analog result for \ref{0z.7} consistently fails, see
\cite{Sh:186} or \cite{Sh:667}; that is: if G.C.H., $\chi > \text{
cf}(\chi) = \kappa$ then we can force a non-reflecting stationary $S
\subseteq S^{\chi^+}_\kappa$ such that the diamond on $S$ fails and
cardinalities and cofinalities are preserved; also G.C.H. continue to
hold.  But
if $\chi$ is strong limit, $\lambda = \chi^+ = 2^\chi$, still we
know something on guessing equalities for every stationary $S \subseteq
S^\lambda_\kappa$; see \cite{Sh:667}.

Note that this $S$ (by \cite{Sh:186}, \cite{Sh:667}) 
in some circumstances has to be ``small"
\mn
\begin{enumerate}
\item[$(*)$]   if ($\chi$ is singular, $2^\chi = \chi^+ =
\lambda,\kappa = \text{ cf}(\chi)$ and) we have the square
$\square_\chi$ (i.e. there exists a sequence $\langle C_\delta:\delta
< \lambda,\delta$ is a limit ordinal$\rangle$, so that $C_\delta$ is
closed and unbounded in $\delta$, cf$(\delta) < \chi \Rightarrow
|C_\delta| < \chi$ and if $\gamma$ is a limit point of $C_\delta$ then
$C_\gamma = C_\delta \cap \gamma$)
\then \, $\diamondsuit_{S^\lambda_\kappa}$ holds, moreover, if $S \subseteq
S^\lambda_\kappa$ reflects in a stationary set of $\delta < \lambda$
then $\diamondsuit_S$ holds, see \cite[\S3]{Sh:186}.
\end{enumerate}
\mn
Also note that our results are of the form ``$\diamondsuit_S$ for
every stationary $S \subseteq S^*$" for suitable $S^* \subseteq
\lambda$.  Now usually this was deduced from the stronger statement
$\diamondsuit^*_S$.  However, the results on $\diamondsuit^*_S$ cannot
be improved, see \ref{2b.7}.

Also if $\chi$ is regular we cannot improve the result to
$\diamondsuit_{S^\lambda_\chi}$, see \cite{Sh:186} or
\cite{Sh:587}, even assuming G.C.H.  
Furthermore the question on $\diamondsuit_{\aleph_2}$ when
$2^{\aleph_1} = \aleph_2 = 2^{\aleph_0}$ was raised.  
Concerning this we show in \ref{2b.4} that
$\diamondsuit_{S^{\aleph_2}_{\aleph_1}}$ may fail (this works in other
cases, too).
\end{problem}

\begin{question}
\label{0z.11}  
Can we deduce any ZFC result on $\lambda$ strongly inaccessible?

By D\v zamonja-Shelah \cite{DjSh:545} we know that failure of 
SNR helps (SNR stands for strong non-reflection), a parallel here is
\ref{ya.4}(2). 

For $\lambda = \lambda^{<\lambda} = 2^\mu$ weakly inaccessible we know
less, still see \cite{Sh:775}, \cite{Sh:898} getting a weaker relative
of diamond, see Definition \ref{2b.9}(2).  Again failure of SNR helps.

On consistency results on SNR see Cummings-D\v zamonja-Shelah
\cite{CDSh:571}, D\v zamonja-Shelah \cite{DjSh:691}.

We thank the audience in the Jerusalem Logic Seminar for their
comments in the lecture on this work in Fall 2006.  
We thank the referee for many helpful remarks.  In particular, as he
urged for details, the proof of \ref{ya.5}(1),(2) just said ``like
\ref{ya.4}", instead we make the old proof of \ref{ya.4}(1) prove
\ref{ya.5}(1),(2) by minor changes and make explicit
\ref{ya.5}(3) which earlier was proved by ``repeat of the second half of the
proof of \ref{ya.4}(1)".  We thank Shimoni Garti for his help in the
proofreading.
\end{question}
 
\begin{notation}
\label{0z.13} 
1) If $\kappa = \text{ cf}(\kappa) < \lambda 
= \text{ cf}(\lambda)$ then we let 
$S^\lambda_\kappa := \{\delta < \lambda:\text{ cf}(\delta) = \kappa\}$.

\noindent
2) ${\cD}_\lambda$ is the club filter on $\lambda$ for $\lambda$ a
   regular uncountable cardinal.
\end{notation}
\newpage

\section {Diamond on successor cardinals} 

Recall (needed only for part (2) of the Theorem \ref{ya.4}):
\begin{definition}
\label{ya.1}  
1) We say $\lambda$ on $S$ has
$\kappa$-SNR or SNR$(\lambda,S,\kappa)$ or $\lambda$ has strong
 non-reflection for $S$ in $\kappa$ \when \, $S \subseteq
 S^\lambda_\kappa : = \{\delta < \lambda:\text{cf}(\delta) = \kappa\}$ 
so $\lambda = \text{ cf}(\lambda) > \kappa = \text{ cf}(\kappa)$ 
and there are $h:\lambda \rightarrow
\kappa$ and a club $E$ of $\lambda$ such that for every $\delta \in
S \cap E$ for some club $C$ of $\delta$, the function $h \restriction C$ is
 one-to-one and even increasing; (note that \wilog \, $\alpha \in \text{
 nacc}(E) \Rightarrow \alpha$ successor and \wilog \, $E = \lambda$,
 so $\mu \in \text{ Reg } \cap \lambda \backslash \kappa^+ \Rightarrow
 \text{ SNR}(\mu,S \cap \mu,\kappa))$.
 If $S = S^\lambda_\kappa$ we may omit it.
\end{definition}

\begin{remark}
\label{ya.2}  Note that by Fodor's lemma if cf$(\delta) = \kappa >
 \aleph_0$ and $h$ is a function from some set $\supseteq \delta$ and
 the range of $h$ is $\subseteq \kappa$ \then \, the following
 conditions are equivalent:
\mn
\begin{enumerate}
\item[$(a)$]   $h$ is one-to-one on some club of $\delta$
\sn
\item[$(b)$]   $h$ is increasing on some club of $\delta$
\sn
\item[$(c)$]   Rang$(h \rest S)$ is unbounded in $\kappa$ for every
 stationary subset $S$ of $\delta$.
\end{enumerate}
\end{remark}

Our main theorem is:
\begin{claim}
\label{ya.4}  Assume $\lambda = 2^\chi = \chi^+$.

\noindent
1) If $S \subseteq \lambda$ is stationary and $\delta \in S
\Rightarrow \text{\rm cf}(\delta) \ne \text{\rm cf}(\chi)$ \then \,
$\diamondsuit_S$ holds.

\noindent
2) If $\aleph_0 < \kappa = \text{\rm cf}(\chi) < \chi$ and $\diamondsuit_S$
fails where $S = S^\lambda_\kappa$ (or just $S \subseteq S^\lambda_\kappa$ is a
stationary subset of $\lambda$) \then \, we have 
{\rm SNR}$(\lambda,\kappa)$ or just $\lambda$ 
has strong non-reflection for $S \subseteq S^\lambda_\kappa$ in $\kappa$.
\end{claim}

\begin{definition}
\label{ya.6} 
1) For a filter $D$ on a set $I$ let 
Dom$(D) := I$ and $S$ is called $D$-positive when $S \subseteq
 I \wedge (I \backslash S) \notin D$ and $D^+ = \{S \subseteq \text{
 Dom}(D):S$ is $D$-positive$\}$ and we let $D+A = \{B \subseteq
   I:B \cup (I \backslash A) \in D\}$ (so if $D = {\cD}_\lambda$,
   the club filter on the regular uncountable $\lambda$ then $D^+$ is
   the family of stationary subsets of $X$).

\noindent
2) For $D$ a filter on a regular uncountable cardinal $\lambda$ which
   extends the club filter, let $\diamondsuit_D$ means: there is $\bar f =
   \langle f_\alpha:\alpha \in S\rangle$ which is a diamond sequence
   for $D$ (or a $D$-diamond sequence) 
which means that $S \in D^+$ and for every $g \in
   {}^\lambda \lambda$ the set $\{\alpha < \lambda:g \rest \alpha =
   f_\alpha\}$ belongs to $D^+$; so $\bar f$ is also a diamond
   sequence for the filter $D+S$, (clearly $\diamondsuit_S$ means
   $\diamondsuit_{{\cD}_\lambda +S}$ for $S$ a stationary subset of
   the regular uncountable $\lambda$).
\end{definition}

\bigskip

\noindent
A somewhat more general version of the theorem is
\begin{claim}
\label{ya.5}  
1) Assume $\lambda = \chi^+ = 2^\chi$ and
$D$ is a $\lambda$-complete filter on $\lambda$ which extends the
club filter.  If $S \in D^+$ and $\delta \in S \Rightarrow 
\text{\rm cf}(\delta) \ne \text{\rm cf}(\chi)$ \then \, we have
$\diamondsuit_{D+S}$.

\noindent
2) We have $\diamondsuit_D$ \when \,:
\mn
\begin{enumerate}
\item[$(a)$]  $\lambda = \lambda^{< \lambda}$
\sn
\item[$(b)$]   $\bar f = \langle f_\alpha:\alpha < \lambda\rangle$ 
lists $\cup\{{}^\alpha \lambda:\alpha < \lambda\}$
\sn
\item[$(c)$]   $S \in D^+$
\sn
\item[$(d)$]   $\bar u = \langle u_\alpha:\alpha \in S\rangle$ and
$u_\alpha \subseteq \alpha$ for every $\alpha \in S$
\sn
\item[$(e)$]  $\chi = \sup\{|u_\alpha|^+:\alpha < \lambda\} < \lambda$
\sn
\item[$(f)$]   $D$ is a $\chi^+$-complete filter on $\lambda$
extending the club filter
\sn
\item[$(g)$]   $(\forall g \in
{}^\lambda \lambda)(\exists^{D^+} \delta \in S)[\delta = \sup\{\alpha
\in u_\delta:g \rest \alpha \in \{f_\beta:\beta \in u_\delta\}\}]$.
\end{enumerate}
\mn
3) Assume $\lambda = \chi^+ = 2^\chi$ and $\aleph_0 < \kappa = 
\text{\rm cf}(\chi) < \chi,S \subseteq S^\lambda_\kappa$ is stationary, and
$D$ is a $\lambda$-complete filter extending the club filter on
$\lambda$ to which $S$ belongs.  If $\diamondsuit_D$ fails
\then \, {\rm SNR}$(\lambda,S,\kappa)$.
\end{claim}

\begin{PROOF}{\ref{ya.4}}
\underline{Proof of \ref{ya.4}}  
Part (1) follows from \ref{ya.5}(1)
for $D$ the filter ${\cD}_\lambda +S$.  Part (2) follows from
\ref{ya.5}(3) for $D$ the filter ${\cD}_\lambda + S$.
\end{PROOF}

\begin{PROOF}{\ref{ya.5}(1)}
\underline{Proof of \ref{ya.5}} 

\noindent
\underline{Proof of part (1)}   

Clearly we can assume
\mn
\begin{enumerate}
\item[$\circledast_0$]  $\chi > \aleph_0$
as for $\chi = \aleph_0$ the statement is empty.
\end{enumerate}
\mn
Let
\mn
\begin{enumerate}
\item[$\circledast_1$]   $\langle f_\alpha:\alpha < \lambda\rangle$
list the set $\{f:f$ is a function from $\beta$ to $\lambda$ for some
 $\beta < \lambda\}$.
\end{enumerate}
\mn
For each $\alpha < \lambda$ clearly $|\alpha| \le \chi$ so let
\mn
\begin{enumerate}
\item[$\circledast_{2,\alpha}$]   $\langle
u_{\alpha,\varepsilon}:\varepsilon < \chi\rangle$ be
$\subseteq$-increasing continuous with union $\alpha$ such that
$\varepsilon < \chi \Rightarrow |u_{\alpha,\varepsilon}| \le \aleph_0 +
|\varepsilon| < \chi$.
\end{enumerate}
\mn
For $g \in {}^\lambda \lambda$ let $h_g \in {}^\lambda \lambda$ be
defined by
\mn
\begin{enumerate}
\item[$\circledast_{3,g}$]  $h_g(\alpha) = \text{ Min}\{\beta <
\lambda:g \restriction \alpha = f_\beta\}$.
\end{enumerate}
\mn
Let cd, $\langle \text{cd}_\varepsilon:\varepsilon < \chi\rangle$ be
such that 
\mn
\begin{enumerate}
\item[$\circledast_4$]  $(a) \quad$ cd is a one-to-one function
from ${}^\chi\lambda$ onto $\lambda$ such that 

\hskip25pt cd$(\bar\alpha) \ge \sup\{\alpha_\varepsilon:\varepsilon <
\chi\}$ (when $\bar\alpha = \langle \alpha_\varepsilon:\varepsilon <
\chi\rangle$) 
\sn
\item[${{}}$]  $(b) \quad$ for $\varepsilon < \chi$,
cd$_\varepsilon$ is a function from $\lambda$ to $\lambda$ such that

\hskip25pt $\bar\alpha = \langle \alpha_\varepsilon:\varepsilon <
\chi\rangle \in {}^\chi\lambda \Rightarrow 
\text{ cd}_\varepsilon(\text{cd}(\bar\alpha)) = \alpha_\varepsilon$
\end{enumerate}
\mn
(they exist as $\lambda = \lambda^\chi$, in the present case
this holds as $2^\chi = \chi^+ = \lambda$).

Now we let (for $\beta < \lambda,\varepsilon < \chi$):
\mn
\begin{enumerate}
\item[$\circledast_5$]   $f^1_{\beta,\varepsilon}$ be the function
from Dom$(f_\beta)$ into $\lambda$ defined by 
$f^1_{\beta,\varepsilon}(\alpha) = \text{ cd}_\varepsilon
(f_\beta(\alpha))$ so Dom$(f^1_{\beta,\varepsilon}) =
\text{ Dom}(f_\beta)$.
\end{enumerate}
\mn
Without loss of generality
\mn
\begin{enumerate}
\item[$\circledast_6$]  $\alpha \in S \Rightarrow \alpha$ is a
limit ordinal.
\end{enumerate}
\mn
For $g \in {}^\lambda \lambda$ and $\varepsilon < \chi$ we let

\begin{equation*}
\begin{array}{clcr}
S^\varepsilon_g = \Big\{ \delta \in S:&\,\delta = \sup\{\alpha \in
u_{\delta,\varepsilon}:\text{ for some } \beta \in u_{\delta,\varepsilon} \\
  &\text{ we have } g \restriction \alpha =f^1_{\beta,\varepsilon}\} \Big\}.
\end{array}
\end{equation*}
\mn
Next we shall show
\mn
\begin{enumerate}
\item[$\circledast_7$]  for some $\varepsilon(*) < \chi$ for every
$g \in {}^\lambda \lambda$ the set $S^{\varepsilon(*)}_g$ is a $D$-positive
 subset of $\lambda$.
\end{enumerate}
\mn
\underline{Proof of $\circledast_7$}   
Assume this fails, so for 
every $\varepsilon < \chi$ there is
$g_\varepsilon \in {}^\lambda \lambda$ such that
$S^\varepsilon_{g_\varepsilon}$ is not $D$-positive and let
$E_\varepsilon$ be a member of $D$ disjoint to
$S^\varepsilon_{g_\varepsilon}$.  Define $g \in {}^\lambda \lambda$ by
$g(\alpha) := \text{ cd}(\langle g_\varepsilon(\alpha):\varepsilon <
\chi\rangle)$ and let $h_g \in {}^\lambda \lambda$ be as in
$\circledast_{3,g}$, i.e. $h_g(\alpha) = \text{ Min}\{\beta:g
\restriction \alpha = f_\beta\}$.

Let $E_* = \{\delta < \lambda:\delta$ is a limit ordinal such that
$\alpha < \delta \Rightarrow h_g(\alpha) < \delta\}$, clearly it is a
club of $\lambda$ hence it belongs to $D$ and so 
$E = \cap\{E_\varepsilon:\varepsilon < \chi\} \cap E_*$ belongs to $D$
as $D$ is $\lambda$-complete and $\chi + 1 < \lambda$.

As $S$ is a $D$-positive subset of $\lambda$ there is $\delta_* \in E
\cap S$.  For each $\alpha < \delta_*$ as $\delta_* \in E \subseteq E_*$
 clearly $h_g(\alpha) <
\delta_*$ and $\alpha$ as well as $h_g(\alpha)$ 
belong to $\cup\{u_{\delta_*,\varepsilon}:\varepsilon < \chi\} = \delta_*$, but
$\langle u_{\delta_*,\varepsilon}:\varepsilon < \chi\rangle$ is
$\subseteq$-increasing hence $\varepsilon_{\delta_*,\alpha} = \text{
min}\{\varepsilon:\alpha \in u_{\delta_*,\varepsilon}$ and 
$h_g(\alpha) \in u_{\delta_*,\varepsilon}\}$
is not just well defined but also $\varepsilon \in
[\varepsilon_{\delta_*,\alpha},\chi) \Rightarrow \{\alpha,h_g(\alpha)\}
\subseteq u_{\delta_*,\varepsilon}$.  
As cf$(\delta_*) \ne \text{ cf}(\chi)$, by an
assumption on $S$, it follows that for some $\varepsilon(*) < \chi$ the
set $B := \{\alpha < \delta_*:\varepsilon_{\delta_*,\alpha} <
\varepsilon(*)\}$ is unbounded below $\delta_*$.

So
\mn
\begin{enumerate}
\item[$(a)$]  $\alpha \in B \Rightarrow \{\alpha,h_g(\alpha)\}
\subseteq u_{\delta_*,\varepsilon(*)}$ and
\sn
\item[$(b)$]   $\alpha \in B \Rightarrow g \rest \alpha =
f_{h_g(\alpha)} \Rightarrow \bigwedge\limits_{\varepsilon < \chi} 
[g_\varepsilon \rest \alpha = f^1_{h_g(\alpha),\varepsilon}] \Rightarrow
g_{\varepsilon(*)} \rest \alpha = f^1_{h_g(\alpha),\varepsilon(*)}$.
\end{enumerate}
\mn
But $\delta_* \in E \subseteq E_{\varepsilon(*)}$ hence $\delta_*
\notin S^{\varepsilon(*)}_{g_{\varepsilon(*)}}$ by the choice of
$E_{\varepsilon(*)}$, but by (a) + (b) and the definition of
$S^{\varepsilon(*)}_{g_{\varepsilon(*)}}$ recalling $\delta_* \in S$
we have sup$(B) = \delta_*
\Rightarrow \delta_* \in S^{\varepsilon(*)}_{g_{\varepsilon(*)}}$,
(where $h_g(\alpha)$ plays the role of $\beta$ in the definition of
$S^\varepsilon_g$ above), contradiction.  So the proof of
$\circledast_7$ is finished.  

Let $\chi_* = (|\varepsilon(*)| + \aleph_0)$ hence $\delta \in S
\Rightarrow |u_{\delta,\varepsilon(*)}|\le \chi_*$ and $\chi^+_* <
\lambda$ as $\chi_* < \chi < \lambda$ because $\aleph_0,\varepsilon(*) <
\chi < \lambda$.  Now we apply \ref{ya.5}(2) which is proved below with
$\lambda,S,D,\chi^+_*,\langle f^1_{\beta,\varepsilon(*)}:\beta <
\lambda\rangle,\langle u_{\delta,\varepsilon(*)}:\delta \in S\rangle$
here standing for $\lambda,S,D,\chi,\bar f,\bar u$ there.  The
conditions there are satisfied hence also the conclusion which says that
$\diamondsuit_D$ holds.
\end{PROOF}

\begin{PROOF}{\ref{ya.5}(2)}
\underline{Proof of \ref{ya.5}(2)}  

Let
\mn
\begin{enumerate}
\item[$\boxtimes_1$]  $\langle \text{cd}_\varepsilon:\varepsilon <
\chi\rangle$ and cd be as in $\circledast_4$ in the proof of part (1),
possible as we are assuming $\chi < \lambda = \lambda^{<\lambda}$
\sn
\item[$\boxtimes_2$]   for $\beta < \lambda$ and $\zeta < \chi$ let
$f^2_{\beta,\zeta}$ be the function with domain Dom$(f_\beta)$ such
that $f^2_{\beta,\zeta}(\alpha) = \text{ cd}_\zeta(f_\beta(\alpha))$
\sn
\item[$\boxtimes_3$]  for $g \in {}^\lambda \lambda$ define $h_g
\in {}^\lambda \lambda$ as in $\circledast_3$ in the proof of part
(1), i.e. $h_g(\alpha) = \text{ Min}\{\beta:g \rest \alpha =
f_\beta\}$. 
\end{enumerate}
\mn
If $2^{< \chi} < \lambda$ our life is easier but we do not
assume this.  For $\delta \in S$ let $\xi^*_\delta$ be a cardinal, and
let $\langle (\alpha^1_{\delta,\xi},\alpha^2_{\delta,\xi}):\xi <
\xi^*_\delta\rangle$ list the set $\{(\alpha_1,\alpha_2) \in
u_\delta \times u_\delta: \text{ Dom}(f_{\alpha_2}) = \alpha_1\}$,
note that $\xi^*_\delta < \chi$, recalling
$|u_\delta| < \chi$ by clause (e) of the assumption.
We now try to choose $(\bar v_\varepsilon,
g_\varepsilon,E_\varepsilon)$ by induction on $\varepsilon < \chi$,
(note that $\bar v_\varepsilon$ is defined from $\langle g_\zeta:\zeta
< \varepsilon\rangle$ (see clause (e) of $\boxtimes_4$ below) so we choose
just $(g_\varepsilon,E_\varepsilon))$, such that:
\mn
\begin{enumerate}
\item[$\boxtimes_4$]   $(a) \quad E_\varepsilon$ is a member of $D$ 
and $\langle E_\zeta:\zeta \le \varepsilon\rangle$ is
$\subseteq$-decreasing with $\zeta$
\sn
\item[${{}}$]   $(b) \quad \bar v_\varepsilon = \langle
v^\varepsilon_\delta:\delta \in S \cap E'_\varepsilon\rangle$ when
$E'_\varepsilon = \cap\{E_\zeta:\zeta < \varepsilon\} \cap \lambda$ so
is $\lambda$ if $\varepsilon = 0$
\sn
\item[${{}}$]   $(c) \quad \langle v^\zeta_\delta:\zeta \le
\varepsilon\rangle$ is $\subseteq$-decreasing with $\zeta$ for each
$\delta \in S \cap E'_\varepsilon$
\sn
\item[${{}}$]   $(d) \quad g_\varepsilon \in{}^\lambda \lambda$
\sn
\item[${{}}$]  $(e) \quad v^\varepsilon_\delta = \{\xi <
\xi^*_\delta$: if $\zeta < \varepsilon$ then $g_\zeta \restriction
\alpha^1_{\delta,\xi} = f^2_{\alpha^2_{\delta,\xi},\zeta}\}$

\hskip25pt (so if $\varepsilon$ is a limit ordinal then $v^\varepsilon_\delta =
\bigcap\limits_{\zeta < \varepsilon} v^\zeta_\delta$ and $\varepsilon = 0
\Rightarrow v^\varepsilon_\delta = \xi^*_\delta$)
\sn
\item[${{}}$]   $(f) \quad$ if $\delta \in E'_\varepsilon \cap
S$ \then \, $v^{\varepsilon +1}_\delta \subsetneqq v^\varepsilon_\delta$ or
$\delta > \sup\{\alpha^1_{\delta,\xi}:\xi \in v^{\varepsilon +1}_\delta\}$.
\end{enumerate}
\mn
Next
\mn
\begin{enumerate}
\item[$\oplus_1$]  we cannot carry the induction, that is for all 
$\varepsilon < \chi$.
\end{enumerate}
\mn
Why?  Assume toward contradiction that $\langle (\bar
v_\varepsilon,g_\varepsilon,E_\varepsilon):\varepsilon < \chi\rangle$
is well defined.
Let $E := \cap\{E_\varepsilon:\varepsilon < \chi\}$, it is a member of
$D$ as $D$ is $\chi^+$-complete. 
Define $g \in {}^\lambda\lambda$ by $g(\alpha) := \text{ cd}
(\langle g_\varepsilon(\alpha):\varepsilon < \chi\rangle)$.  
Let $E_* = \{\delta < \lambda:\delta$ a limit ordinal such that
$h_g(\alpha) < \delta$ and $\delta > \sup(\text{Dom}(f_\alpha) \cup
\text{ Rang}(f_\alpha))$ for every $\alpha < \delta\}$, so $E_*$ is a club of
$\lambda$ hence it belongs to $D$.  By assumption (g) of the claim the set

\[
S_g := \{\delta \in S:\delta = \sup\{\alpha \in u_\delta:(\exists \beta \in
u_\delta)(f_\beta = g \restriction \alpha)\}\}
\]
\mn
is $D$-positive, so we can choose $\delta \in E \cap E_* \cap
S_g$.  Hence $B := \{\alpha \in u_\delta:(\exists \beta \in u_\delta)
(f_\beta = g \rest \alpha)\}$ is an unbounded subset
of $u_\delta$ and let $h:B \rightarrow u_\delta$ be 
$h(\alpha) = \text{ min}\{\beta \in u_\delta:f_\beta = g \rest
\alpha\}$, clearly $h$ is a function from $B$ into
$u_\delta$.  Now $\alpha \in B \wedge \zeta < \chi \Rightarrow
f_{h(\alpha)} = g \rest \alpha \wedge \zeta < \chi \Rightarrow 
f^2_{h(\alpha),\zeta} = g_\zeta \rest \alpha$,
so for $\alpha \in B$ the pair $(\alpha,h(\alpha))$ belongs to
$\{(\alpha^1_{\delta,\xi},\alpha^2_{\delta,\xi}):\xi \in
v^\varepsilon_\delta\}$ for every $\varepsilon < \chi$.  Hence for any
$\varepsilon < \chi$ we have $B
\subseteq \{\alpha^1_{\delta,\xi}:\xi \in v^\varepsilon_\delta\}$ so
$\delta = \sup\{\alpha^1_{\delta,\xi}:\xi \in v^\varepsilon_\delta\}$.

So for the present $\delta$, in clause (f) of $\boxtimes_4$ the second
possibility never occurs. 

So clearly $\langle v^\varepsilon_{\delta_*}:
\varepsilon < \chi\rangle$ is strictly $\subseteq$-decreasing, i.e. is
$\subset$-decreasing which is impossible as 
$|v^0_{\delta_*}| = \xi^*_{\delta_*} < \chi$.  So
we have proved $\oplus_1$ hence we can assume
\mn
\begin{enumerate}
\item[$\oplus_2$]  there is $\varepsilon < \chi$ such that
we have defined our triple for every $\zeta < \varepsilon$ but we
cannot define for $\varepsilon$.
So we have $\langle(\bar v_\zeta,g_\zeta,E_\zeta):\zeta <
\varepsilon\rangle$.
\end{enumerate}
\mn
As in $\boxplus_4(e)$, let
\mn
\begin{enumerate}
\item[$\odot_1$]   $E'_\varepsilon$ be $\lambda$ if $\varepsilon =
0$ and $\cap\{E_\zeta:\zeta < \varepsilon\}$ if $\varepsilon > 0$ and let
$S_* := S \cap E'_\varepsilon$.
\end{enumerate}
\mn
Clearly $\bar v_\varepsilon$ is well defined, see clauses (b),(e) of
$\boxtimes_4$, and for $\delta \in S_*$
let ${\cF}_\delta = \{f^2_{\alpha^2_{\delta,\xi},\varepsilon}:\xi
\in v^\varepsilon_\delta\}$, so each member is a function from some
$\alpha \in u_\delta \subseteq \delta$ into some ordinal $< \delta$.

Let
\mn
\begin{enumerate}
\item[$\odot_2$]   $S^*_1 := \{\delta \in S_*:
\text{ there are } f',f'' \in {\cF}_\delta
\text{ which are incompatible as functions}\}$
\sn
\item[$\odot_3$]   $S^*_2 := \{\delta \in S_*:
\delta \notin S^*_1 \text{ but the function }
\cup\{f:f \in {\cF}_\delta\} \text{ has domain } \ne \delta\}$
\sn
\item[$\odot_4$]   $S^*_3 = S_* \backslash(S^*_1 \cup S^*_2)$.
\end{enumerate}
\mn
For $\delta \in S^*_3$ let $g^*_\delta = \cup\{f:f \in 
{\cF}_\delta\}$, so by the definition of $\langle
S^*_\ell:\ell=1,2,3\rangle$ clearly $g^*_\delta \in {}^\delta \delta$.
Now if $\langle g^*_\delta:\delta \in S^*_3\rangle$ is a diamond
sequence for $D$ then we are done.

So assume that this fails, so for some $g \in {}^\lambda\lambda$ and member
 $E$ of $D$ we have $\delta \in S^*_3 \cap E \Rightarrow
g^*_\delta \ne g \restriction \delta$.  Without loss of generality $E$
is included in $E'_\varepsilon$.  But then we could have chosen
$(g,E)$ as $(g_\varepsilon,E_\varepsilon)$, recalling 
$\bar v_\varepsilon$ was already chosen.
Easily the triple $(g_\varepsilon,E_\varepsilon,\bar v_\varepsilon)$
is as required in $\oplus_1$, contradicting the
choice of $\varepsilon$ in $\oplus_2$ 
so we are done proving part (2) of Theorem \ref{ya.4} hence also part (1).
\end{PROOF}

\begin{proof}
\underline{Proof of part (3)}   

We use cd,cd$_\varepsilon$ (for
$\varepsilon < \chi$), $\left< \langle u_{\alpha,\varepsilon}:\varepsilon <
\chi\rangle:\alpha < \lambda \right>,\langle f_\alpha:\alpha <
\lambda\rangle,\langle f^1_{\alpha,\varepsilon}:\alpha <
\lambda,\varepsilon < \chi\rangle$ and $S^\varepsilon_g$ for
$\varepsilon < \kappa$ as in the proof of part (1).

Recall $\kappa$, a regular uncountable cardinal, is the cofinality of the
singular cardinal $\chi$ and let $\langle
\chi_\gamma:\gamma < \kappa\rangle$ be increasing with limit $\chi$.
For every $\gamma < \kappa$ we ask:
\mn
\underline{The $\gamma$-Question}:  
Do we have: for every $g \in {}^\lambda \lambda$, the
following is a $D$-positive subset of $\lambda$:

\noindent
$\{\delta \in S:S_\gamma[g] \cap \delta$ is a stationary subset of
$\delta\}$ where $S_\gamma[g] := \{\zeta < \lambda$: cf$(\zeta)
\in [\aleph_0,\kappa)$, sup$(u_{\zeta,\chi_\gamma}) = \zeta$
and for arbitrarily large $\alpha \in u_{\zeta,\chi_\gamma}$ 
for some $\beta \in u_{\zeta,\chi_\gamma}$ and
$\varepsilon < \chi_\gamma$, we have Dom$(f_\beta) = \alpha$ and $g
\restriction \alpha = f^1_{\beta,\varepsilon}\}$.
\bigskip

\noindent
\underline{Case 1}:  For some $\gamma < \kappa$, the answer is yes.

Choose $\langle C_\delta:\delta \in S\rangle$ such that $C_\delta$ is
a club of $\delta$ of order type cf$(\delta) = \kappa$.

For $\delta \in S \subseteq S^\lambda_\kappa$ let $u_\delta :=
\cup\{u_{\alpha,\chi_\gamma}:\alpha \in C_\delta\}$.

Clearly
\mn
\begin{enumerate}
\item[$\boxplus_2$]   $|u_\delta| \le \kappa + \chi_\gamma < \chi$
\sn
\item[$\boxplus_3$]    for every $g \in {}^\lambda \lambda$ 
for $D$-positively many $\delta \in S$, we have
$\delta = \sup\{\alpha \in u_\delta:g \rest \alpha \in
\{f^1_{\beta,\varepsilon}:\varepsilon < \chi_\gamma$ and 
$\beta \in u_\delta\}\}$.
\end{enumerate}
\mn
Why $\boxplus_3$ holds?  Given $g \in {}^\lambda \lambda$, let $h_g
\in {}^\lambda \lambda$ be defined by $h_g(\alpha) = \text{ min}
\{\beta < \lambda:g \rest \alpha = f_\beta\}$, so $h_g(\alpha) \ge
\alpha$ (but is less than $\lambda$).  Let $E_g = \{\delta < \lambda:\delta$
is a limit ordinal such that $(\forall \alpha < \delta)h_g(\alpha) <
\delta\}$, so $E_g$ is a club of $\lambda$ and let $E'_g$ be the set of
accumulation points of $E_g$, so $E'_g$, too, is a club of $\lambda$.
By the assumption of this case, the set $S' := \{\delta \in S:\delta
\cap S_\gamma[g]$ is a stationary subset of $\lambda\}$ is $D$-positive,
hence $S'' := S' \cap E'_g$ is a $D$-positive subset of $\lambda$.  Let
$\delta \in S''$, by $E'_g$'s definition, we can find $B^0_\delta
\subseteq E_g \cap \delta$ unbounded in $\delta$, so \wilog \,
$B^0_\delta$ is closed.  But $S_\gamma[g] \cap \delta$
is a stationary subset of $\delta$, recalling $\delta \in S''$, so
$B^1_\delta = B^0_\delta \cap S_\gamma[g] \cap C_\delta$ is a
stationary subset of $\delta$ as $B^0_\delta,C_\delta$ are closed
unbounded subsets of $\delta$. 

Clearly $\zeta \in B^1_\delta \Rightarrow \zeta \in C_\delta
\Rightarrow u_{\zeta,\chi_\gamma} \subseteq u_\delta$ by the
definitions of $B^1_\delta$ and $u_\delta$.  Also $\zeta \in B^1_\delta
\Rightarrow \zeta \in S_\gamma[g] \Rightarrow (\zeta$ is a limit
ordinal) $\wedge \zeta = \sup(u_{\zeta,\chi_\gamma}) = \sup\{\alpha \in
u_{\zeta,\chi_\gamma}:(\exists \beta \in u_{\zeta,\chi_\gamma})
(\exists \varepsilon < \chi_\gamma)(g \rest \alpha =
f^1_{\beta,\varepsilon})\} \Rightarrow$ (($\zeta$ is a limit ordinal)
$\wedge \zeta = \sup\{\alpha \in u_\delta \cap \zeta:(\exists \beta
\in u_\delta \backslash \alpha)(\exists \varepsilon < \chi_\gamma)(g
\rest \alpha = f^1_{\beta,\varepsilon})\})$.

As $B^1_\delta$ is unbounded in $\delta$ being stationary we are
done proving $\boxplus_3$.
\medskip

Now \wilog \, every $\delta \in S$ is divisible by $\chi$ hence
$\delta = \chi_\gamma \delta$ and
let $u'_\delta = u_\delta \cup\{\chi_\gamma \alpha +
\varepsilon:\alpha \in u_\delta,\varepsilon < \chi_\gamma\}$, so
$u_\delta$ is an unbounded subset of $\delta$, and let $f'_\beta =
f^1_{\alpha,\varepsilon}$ when $\beta = \chi_\gamma \alpha
+\varepsilon,\varepsilon < \chi_\gamma$.  So translating what we have is:
\mn
\begin{enumerate}
\item[$\boxplus_4$]   $(a) \quad \langle f'_\alpha:\alpha <
\lambda\rangle$ is a sequence of members of $\cup\{{}^\beta
\lambda:\beta < \lambda\}$
\sn
\item[${{}}$]   $(b) \quad$ for $\delta \in S,u'_\delta$ is 
an unbounded subset of $\delta$ of cardinality 

\hskip25pt $\le \chi_\gamma \times \chi_\gamma = \chi_\gamma (< \chi)$
\sn
\item[${{}}$]   $(c) \quad$ for every $g \in {}^\lambda \lambda$ for
$D$-positively many $\delta \in S$ we have

\hskip25pt  $\delta = \sup\{\alpha
\in u'_\delta:(\exists \beta \in u'_\delta)(g \rest \alpha = f'_\beta)\}$.
\end{enumerate}
\mn
Now we can apply part (2) with $\langle f'_\alpha:\alpha <
\lambda\rangle,\langle u'_\delta:\delta \in S\rangle$ replacing $\bar
f,\langle u_\delta:\delta \in S\rangle$.

So as there we can prove $\diamondsuit_S$, hence we are done.
\bigskip

\noindent
\underline{Case 2}:  For every $\gamma < \kappa$ the answer is no.

Let $(g_\gamma,E_\gamma)$ exemplify that the answer for $\gamma$ is
no; so $g_\gamma \in {}^\lambda \lambda$ and $E_\gamma \in D$.  
Let $E = \bigcap\limits_{\gamma < \kappa} E_\gamma$, so $E$ is a member of
$D$.   Let $g \in {}^\lambda \lambda$ be defined by $g(\alpha) =
\text{ cd}(\langle g_\gamma(\alpha):\gamma < \kappa\rangle \char 94
(0)_\chi)$, i.e. cd$_\varepsilon(g(\alpha))$ is $g_\gamma(\alpha)$ if
$\gamma < \kappa$ and is $0$ if $\varepsilon \in [\kappa,\chi)$.

Let

\begin{equation*}
\begin{array}{clcr}
E_g := \{\delta < \lambda:&\,\delta \text{ a limit ordinal such that if }
\alpha < \lambda \text{ then } h_g(\alpha) < \delta \\
  &\text{and } \delta > \sup(\text{Dom}(f_\alpha) \cup \text{
  Rang}(f_\alpha))\}.
\end{array}
\end{equation*}
\mn
We now define $h:\lambda \rightarrow \kappa$ as follows
\mn
\begin{enumerate}
\item[$\boxplus_5$]  for $\beta < \lambda$
\begin{enumerate}
\item[$(a)$]   if cf$(\beta) \notin [\aleph_0,\kappa)$ or $\beta
\notin E_g$ then $h(\beta) = 0$
\sn
\item[$(b)$]  otherwise
\end{enumerate}

\begin{equation*}
\begin{array}{clcr}
h(\beta) = \text{ min}\{\gamma < \kappa:&\beta = \sup\{\alpha_1 \in
u_{\beta,\chi_\gamma}:\text{ for some } \alpha_2 \in u_{\beta,\chi_\gamma} \\
  &\text{ and } \varepsilon < \chi_\gamma \text{ we have } g
\restriction \alpha_1 = f^1_{\alpha_2,\varepsilon}\}\}.
\end{array}
\end{equation*}
\end{enumerate}
\mn
Now
\mn
\begin{enumerate}
\item[$\boxplus_6$]   $h:\lambda \rightarrow \kappa$ is well defined.
\end{enumerate}
\mn
Why $\boxplus_6$ holds?  Let $\beta < \lambda$.  If cf$(\beta) \notin 
[\aleph_0,\kappa)$ or
$\beta \notin E_g$ then $h(\alpha) = 0 < \kappa$ by clause (a) of
$\boxplus_5$.  So assume cf$(\beta) \in [\aleph_0,\kappa)$ and $\beta
\in E_g$.  Let $\langle \gamma^1_{\beta,\varepsilon}:\varepsilon <
\text{ cf}(\beta)\rangle$ be increasing with limit $\beta$ and let
$\gamma^2_{\beta,\varepsilon} = \text{ min}\{\gamma:g \rest
\gamma^1_{\beta,\varepsilon} = f_\gamma\}$, so $\varepsilon < \text{
cf}(\beta) \Rightarrow \gamma^2_{\beta,\varepsilon} < \beta$ as $\beta
\in E_g$.  But $\langle u_{\beta,\chi_\zeta}:\zeta < \text{
cf}(\chi)\rangle$ is $\subseteq$-increasing with union $\beta$ so for each
$\varepsilon < \text{ cf}(\beta)$ there is $\zeta =
\zeta_{\beta,\varepsilon} < \text{ cf}(\chi)$ such that
$\{\gamma^1_{\beta,\varepsilon},\gamma^2_{\beta,\varepsilon}\}
\subseteq u_{\beta,\chi_\zeta}$.  As cf$(\beta) < \kappa = \text{
cf}(\chi)$ for some $\zeta < \kappa$ the set $\{\varepsilon < \text{
cf}(\beta):\zeta_{\beta,\varepsilon} < \zeta\}$ is unbounded in
cf$(\beta)$.  So $\zeta$ can serve as $\gamma$ in clause (b) of
$\boxplus_5$ so $h(\beta)$ is well defined, in particular is less than
$\kappa$ so we have proved $\boxplus_6$.
\mn
\begin{enumerate}
\item[$\boxplus_7$]   if $\delta \in S \cap E_\gamma$ then for
some club $C$ of $\delta$ the function $h \restriction C$ is increasing.
\end{enumerate}
\mn
Why $\boxplus_7$ holds?  If not, then by Fodor's 
lemma for some $\gamma < \kappa$ the
set $\{\delta' \in \delta \cap S:h(\delta') \le \gamma\}$ is a
stationary subset of $\delta$, and we get contradiction to the choice
of $E_\gamma$ so $\boxplus_7$ holds indeed.

So $h$ is as promised in the claim.  
\end{proof}

Note
\begin{observation}
\label{ya.12}  
If $\kappa_* < \lambda$ are regular, $S^\lambda_{\kappa_*}$ strongly does not
reflect in $\lambda$ for every $\kappa \in \text{\rm Reg } \cap \kappa_*$
and $\Pi(\text{Reg } \cap \kappa_*) < \lambda$, \then \, :
\mn
\begin{enumerate}
\item[$(a)$]   $S^\lambda_{< \kappa_*}$ can be divided to 
$\le \Pi(\text{\rm Reg } \cap \kappa_*)$ sets, each not reflecting
   in any $\delta \in S^\lambda_{< \kappa_*}$

in particular
\sn
\item[$(b)$]   $S^\lambda_{\aleph_0}$ can be divided to $\le
\Pi(\text{\rm Reg } \cap \kappa_*)$ sets each not reflecting in any
$\delta \in S^\lambda_{< \kappa_*}$.
\end{enumerate}
\end{observation}

\begin{remark}
\label{ya.14}  
1) Of course if $\lambda$ has $\kappa$-SNR then this
holds for every regular $\lambda' \in (\kappa,\lambda)$.

\noindent
2) We may state the results, using $\lambda^*_\kappa$ (see below).
\end{remark}

\begin{definition}
\label{ya.16} 
For each regular $\kappa$ let $\lambda^*_\kappa = 
\text{ Min}\{\lambda:\lambda$ regular fails to
 have $\kappa -$ SNR$\}$, and let $\lambda^*_\kappa$ be $\infty$ (or
 not defined) if there is no such $\lambda$.
\end{definition}
\newpage

\section {Consistent failure on $S^2_1$} 

A known question was:
\begin{question}
\label{2b.1}  
For $\theta \in \{\aleph_0,\aleph_1\}$ do
we have $(2^{\aleph_0} = 2^{\aleph_1} = \aleph_2 \Rightarrow 
\diamondsuit_{S^{\aleph_2}_\theta})$?

So for $\theta = \aleph_0$ the answer is yes (by \ref{ya.4}(1)), but
what about $\theta = \aleph_1$?  We noted some years ago that easily:
\end{question}

\begin{claim}
\label{2b.4}  Assume $\bold V \models$ {\rm GCH} or even just
$2^{\aleph_\ell} = \aleph_{\ell +1}$ for $\ell = 0,1,2$.  \Then \, some
forcing notion $\Bbb P$ satisfies
\mn
\begin{enumerate}
\item[$(a)$]   $\Bbb P$ is of cardinality $\aleph_3$
\sn
\item[$(b)$]   forcing with $\bbP$ preserves cardinals and
cofinalities  
\sn
\item[$(c)$]  in $\bold V^{\bbP},2^{\aleph_0} = 2^{\aleph_1} = 
\aleph_2,2^{\aleph_2} = \aleph_3$
\sn
\item[$(d)$]   in $\bold V^{\bbP},
\diamondsuit_S$ fails where $S = \{\delta
< \aleph_2:{\text{\rm cf\/}}(\delta) = \aleph_1\}$, moreover
\begin{enumerate}
\item[$(*)$]   there is a sequence $\bar A =  
\langle A_\delta:\delta \in S\rangle$ where 
$A_\delta$ an unbounded subset of $\delta$ of order type
$\omega_1$ satisfying
\sn
\item[$(**)$]   if $\bar f = \langle f_\delta:\delta \in 
S\rangle,f_\delta \in {}^{(A_\delta)}(\omega_1)$, \then \, there is $f
\in {}^{(\omega_2)}(\omega_1)$ such that $\delta \in
S \Rightarrow \delta > \sup(\{\alpha \in A_\delta:
f(\alpha) \le f_\delta(\alpha)\})$.
\end{enumerate}
\end{enumerate}
\end{claim}

\begin{remark}
Similarly for other cardinals.
\end{remark}

\begin{PROOF}{\ref{2b.4}}
There is an $\aleph_1$-complete $\aleph_3$-c.c. forcing
notion $\bbP$ not collapsing cardinals, not changing cofinalities,
preserving $2^{\aleph_\ell} = \aleph_\ell$ for $\ell=0,1,2$ and
$|\bbP| = \aleph_3$ such that in $\bold V^{\bbP}$, we have $(*)$,
in fact more\footnote{I.e. there is $\bar A = \langle A_\delta:\delta \in
S\rangle$ where $A_\delta$ is an unbounded subset of $\delta$ of order
type $\omega_1$ satisfying:
\mn
\begin{enumerate}
\item[$\oplus$]   if $\bar f = \langle f_\delta:\delta \in
S\rangle,f_\delta \in {}^{(A_\delta)}\omega_1$ then there is $f \in
{}^{(\omega_2)}\omega_1$ such that for every $\delta \in
S^{\aleph_2}_{\aleph_1}$ for every $\alpha \in A_\delta$
large enough we have $f(\alpha) = f_\delta(\alpha)$.
\end{enumerate}}
 than $(*)$ holds - see \cite{Sh:587}.  Let $\bbQ$ be the forcing of
adding $\aleph_2$ Cohen or just any 
c.c.c. forcing notion of cardinality $\aleph_2$ adding $\aleph_2$
reals (can be $\name{\bbQ}$, a $\bbP$-name).
\mn
Now we shall show that $\bbP * \bbQ$, equivalently $\bbP \times \bbQ$ 
is as required: 
\medskip

\noindent
\underline{Clause (a)}:  

$|\bbP * \bbQ| = \aleph_3$; trivial.
\smallskip

\noindent
\underline{Clause (b)}:  

Preserving cardinals and cofinalities; obvious as both $\bbP$ and
$\bbQ$ do this.
\smallskip

\noindent
\underline{Clause (c)}:  Easy. 
\smallskip

\noindent
\underline{Clause (d)}:  In $\bold V^{\bbP}$ we have
$(*)$ as exemplified by say $\bar A = \langle A_\delta:\delta \in 
S\rangle$.  We shall show that $\bold V^{\bbP * \bbQ} 
\models ``\bar A$ satisfies $(**)"$.  Otherwise
in $\bold V^{{\bbP}* {\bbQ}}$ we have $\bar f =
\langle f_\delta:\delta \in S\rangle$ say in 
$\bold V[G_{\bbP},G_{\bbQ}]$ a counterexample \then \,
in $\bold V[G_{\bbP}]$ for some $q \in \bbQ$ and $\name{\bar f}$
we have

\begin{equation*}
\begin{array}{clcr}
\bold V[G_{\bbP}] \models (q \Vdash_{\bbQ} 
``\name{\bar f} = &\langle \name f_\delta:\delta \in S\rangle \text{ where } 
{\name f_\delta}:A_\delta \rightarrow \omega_1 \text{ for each } 
\delta \in S \\
  &\text{ form a counterexample to } (*)").
\end{array}
\end{equation*}
\mn
Now in $\bold V[G_{\bbP}]$ we can define 
$\bar g = \langle g^1_\delta:\delta \in S\rangle \in 
\bold V[G_{\bbP}]$ where $g^1_\delta$ a function 
with domain $A_\delta$, by

\[
g^1_\delta(\alpha) = \{i:q \nVdash {\name f_\delta}(\alpha) \ne i\}.
\]
\mn
So in $\bold V[G_{\bbP}]$ we have
$q \Vdash_{\bbQ} ``\bigwedge\limits_{\delta \in S} (\forall \alpha
\in A_\delta){\name f_\delta}(\alpha) \in
g^1_\delta(\alpha)\}"$.  Also $g^1_\delta(\alpha)$ is a countable
subset of $\omega_1$ as $\bbQ$ satisfies the c.c.c.

For $\delta \in S$ we define a function $g_\delta:A_\delta \rightarrow
\omega_1$ by letting $g_\delta(\alpha) = (\sup(g^1_\delta(\alpha)) +1$ hence
$g_\delta(\alpha) < \omega_1$ so $\langle g_\delta:\delta \in
S\rangle$ is as required on $\bar f$ in $(**)$ in 
$\bold V[G_{\bbP}]$, of course.  Apply clause $(**)$ in
$\bold V[G_{\bbP}]$ to $\langle g_\delta:\delta \in S\rangle$
so we can find $g:\omega_2 \rightarrow \omega_1$ such that 
$\bigwedge\limits_{\delta \in S} \, \delta > \sup\{\alpha \in A_\delta,
g_\delta(\alpha) > g(\alpha)\}$.  Now $g$ is as required also in
$\bold V[G_{\bbP}][G_{\bbQ}]$.  
\end{PROOF}

We may wonder can we strengthen the conclusion of
\ref{ya.4} to $\diamondsuit^*_S$ (of
course the demand in clause (e) and (f) in claim \ref{2b.7} below  
are necessary, i.e. otherwise $\diamondsuit^*_S$ holds).  
The answer is \underline{not} as: (the restriction in (e)
and in (f) are best possible).
\begin{observation}
\label{2b.7}  
Assume $\lambda = \lambda^{< \lambda},S \subseteq S^\lambda_\kappa$.

Then for some $\bbP$
\mn
\begin{enumerate}
\item[$(a)$]   $\bbP$ is a forcing notion 
\sn
\item[$(b)$]  $\bbP$ is of cardinality $\lambda^+$ satisfying the
$\lambda^+$-c.c.
\sn
\item[$(c)$]   forcing with $\bbP$ does not collapse cardinals and
does not change cofinality
\sn
\item[$(d)$]   forcing with $\bbP$ adds 
no new $\eta \in {}^{\lambda >}\text{\rm Ord}$
\sn
\item[$(e)$]   $\diamondsuit^*_S$ fails for every stationary subset
$S$ of $\lambda$ such that
\begin{enumerate}
\item[$(\alpha)$]   $S \subseteq S^\lambda_\kappa$ when $(\exists
\mu < \lambda)[\mu^{<\kappa>_{\text{\rm tr}}} = \lambda]$ 

or just
\sn
\item[$(\beta)$]   $\alpha \in S \Rightarrow
|\alpha|^{<\text{\rm cf}(\alpha)>_{\text{\rm tr}}} > |\alpha|$
\end{enumerate}
\item[$(f)$]   $(D \ell)_S$, see below, fails for every $S \subseteq
S^\lambda_\kappa$ when $\alpha \in S \Rightarrow
|\alpha|^{<\text{\rm cf}(\alpha)>_{\text{tr}}} = \lambda$.
\end{enumerate}
\end{observation}

Recalling
\begin{definition}
\label{2b.9}  
1) For $\mu \ge \kappa = \text{ cf}(\kappa)$ let 
$\mu^{<\kappa>_{\text{tr}}} = \{|{\cT}|:{\cT}
\subseteq {}^{\kappa \ge}\mu$ is closed under initial segments (i.e. a
subtree) such that $|{\cT} \cap {}^{\kappa >}\mu| \le \mu\}$.

\noindent
2) For $\lambda$ regular uncountable and stationary $S \subseteq
\lambda$ let $(D \ell)_S$ mean that there is a sequence $\bar{\cP} 
= \langle {\cP}_\delta:\delta \in S\rangle$ witnessing it
 which means: 
\mn
\begin{enumerate}
\item[$(*)_{\bar{\cP}}$]   $(a) \quad {\cP}_\delta \subseteq
 {}^\delta \delta$ has cardinality $< \lambda$
\sn
\item[${{}}$]   $(b) \quad$ for every $f \in {}^\lambda \lambda$ the
set $\{\delta \in S:f \rest \delta \in {\cP}_\delta\}$ is stationary
\end{enumerate}
\mn
(for $\lambda$ successor it is equivalent to $\diamondsuit_S$; for
$\lambda$ strong inaccessible it is trivial).
\end{definition}

\begin{PROOF}{\ref{2b.7}}
\underline{Proof of \ref{2b.7}}  
Use $\bbP =$ adding $\lambda^+,\lambda$-Cohen subsets.

The proof is straight.  
\end{PROOF}
\bigskip

\begin{remark}
\label{2b.11}  The consistency results in \ref{2b.7}
are best possible, see \cite{Sh:829}.
\end{remark}


\end{document}